\documentclass{amsart}
\usepackage{amssymb,amsmath}
\headsep=1cm \numberwithin{equation}{section}
\newtheorem{thm}{Theorem}[section]
\newtheorem{cor}[thm]{Corollary}
\newtheorem{lem}[thm]{Lemma}
\newtheorem{prop}[thm]{Proposition}

\newtheorem{exam}[thm]{Example}

\newcommand{\Ann}{\mbox{Ann}\,}

\newcommand{\Spec}{\mbox{Spec}\,}

\newcommand{\Ass}{\mbox{Ass}\,}
\newcommand{\Assh}{\mbox{Assh}\,}
\newcommand{\Att}{\mbox{Att}\,}
\newcommand{\Supp}{\mbox{Supp}\,}

\renewcommand{\dim}{\mbox{dim}\,}

\newcommand{\h}{\mbox{ht}\,}

\newcommand{\fa}{\mathfrak{a}}
\newcommand{\fb}{\mathfrak{b}}

\newcommand{\fm}{\mathfrak{m}}
\newcommand{\fp}{\mathfrak{p}}
\newcommand{\fq}{\mathfrak{q}}

\newcommand{\LH}{Lichtenbaum-Hartshorne Theorem}
\newcommand{\lc}{\mbox{H}\, }

\begin{document}

\bibliographystyle{amsplain}

\author{Mohammad T. Dibaei}
\address{ Mohammad T. Dibaei\\Faculty of Mathematical Sciences, Teacher Training
University, Tehran, Iran, and Institute for Theoretical Physics and
Mathematics (IPM), Tehran, Iran.}

\email{dibaeimt@ipm.ir}

\author{Raheleh Jafari}
\address{Raheleh Jafari\\Faculty of Mathematical Sciences, Teacher Training
University, Tehran, Iran}

\email{jafarirahele@yahoo.com}

\keywords{Local cohomology, Associated primes, Attached primes\\
The research of the first author was in part supported from IPM (No.
00000000).}

\subjclass[2000]{13D45, 13D07}

\title[Top local cohomology modules with specified attached primes ]
{Top local cohomology modules  \\ with specified attached primes}


\begin{abstract}
 Let $(R,\fm)$ be a complete Noetherian local ring and let $M$ be
a finite $R$--module of positive Krull dimension $n$. It is shown
that any subset $T$ of $\Assh_R(M)$ can be expressed as the set of
attached primes of the top local cohomology module $\lc^n_\fa(M)$
for some ideal $\fa$ of $R$. Moreover if $\fa$ is an ideal of $R$
such that the set of attached primes of $\lc^n_\fa(M)$ is a
non--empty proper subset of $\Assh_R(M)$, then
$\lc^n_\fa(M)\cong\lc^n_\fb(M)$ for some ideal $\fb$ of $R$ with
$\dim_R (R/\fb)=1$.
\end{abstract}

\maketitle

\section{Introduction}

Throughout $(R,\fm)$ is a  commutative Noetherian local ring with
maximal ideal $\fm$, $M$ is a non-zero finite (i.e. finitely
generated) $R$--module with positive Krull dimension $n:=\dim_R(M)$
and $\fa$ denotes an ideal of $R$. Recall that for an $R$--module
$N$, a prime ideal $\fp$ of $R$ is said to be an {\it attached
prime}
 of
$N$, if $\fp=\Ann_R(N/K)$ for some submodule $K$ of $N$ (see
\cite{MS}). The set of attached primes of $N$ is denoted by
$\Att_R(N)$. If $N$ is an Artinian $R$--module so that $N$ admits a
reduced secondary representation $N=N_1+\cdots+N_r$ such that $N_i$
is $\fp_i$--secondary, $i=1,\ldots,r$, then
$\Att_R(N)=\{\fp_1,\ldots,\fp_r\}$ is a finite set.

Denote by $\lc^n_\fa(M)$ the $n$th right derived functor of
$$\Gamma_\fa(M)=\{x\in M|\, \fa^rx=0 \ \mbox{for some positive
integer} \  r \}$$   applied to $M$. It is well-known that
$\lc^n_\fa(M)$ is an Artinian module. Macdonald and Sharp, in
\cite{MS}, studied $\lc^n_\fm(M)$ and showed that
$\Att_R(\lc^n_\fm(M))= \Assh_R(M)$ where $\Assh_R(M):=\{\fp \in
\Ass_R(M)|\, \dim_R(R/\fp)=n\}$. It is shown in \cite[Theorem
A]{DY1}, that for any arbitrary ideal $\fa$ of $R$,
$\Att_R(\lc^n_\fa(M))=\{\fp\in\Ass_R(M)|\, \lc^n_\fa(R/\fp)\neq 0\}$
which is a subset of $\Assh_R(M)$. In \cite{DY2}, the structure of
$\lc^n_\fa(M)$ is studied by the first author and Yassemi and they
showed that, in case $R$ is complete, for any pair of ideals $\fa$
and $\fb$ of $R$, if $\Att_R(\lc^n_\fa(M))=\Att_R(\lc^n_\fb(M))$,
then $\lc^n_\fa(M) \cong \lc^n_\fb(M)$. They also raised the
following question in \cite[Question 2.9]{DY3} which is the main
object of this paper.

\noindent{\bf Question.} For any subset $T$ of $\Assh_R(M)$, is
there an ideal $\fa$ of $R$ such that $\Att_R(\lc^n_\fa(M)) = T$?

This paper provides a positive answer for this question in the case
$R$ is complete.


\section{Main Result}
In this section we assume that $R$ is complete with respect to the
$\fm$--adic topology. As mentioned above, $\Att_R(\lc^n_\fm(M)) =
\Assh_R(M)$ and $\Att_R(\lc^n_R(M)) = \emptyset$ is the empty set.
Also  $\Att_R(\lc^n_\fa(M)) \subseteq \Assh_R(M)$ for all ideals
$\fa$ of $R$. Our aim is to show that as $\fa$ varies over ideals of
$R$, the set $\Att_R(\lc^n_\fa(M))$ takes all possible subsets of
$\Assh_R(M)$ (see Theorem 2.8). In the following results we always
assume that $T$ is a non--empty proper subset of $\Assh_R(M)$

In our first result we find a characterization for a subset of
$\Assh_R (M)$ to be the set of attached primes of the top local
cohomology of $M$ with respect to an ideal $\fa$.

\begin{prop}
Assume that $n:=\dim_R(M)\geq 1$ and that $T$ is a proper non-empty
subset of $\Assh_R(M)$. Set $\Assh_R(M)\setminus
T=\{\fq_1,\ldots,\fq_r\}$. The following statements are equivalent.
\begin{enumerate}
  \item[(i)] There exists an ideal $\fa$ of $R$ such that
  $\Att_R(\lc^n_\fa(M))=T$.
  \item[(ii)] For each $i,\,1\leq i\leq r$, there exists $Q_i\in \Supp_R(M)$
   with $\dim_R(R/Q_i)=1$
  such that
  $$\underset{\fp\in T}{\bigcap}\fp\nsubseteq Q_i \quad \mbox{and}
  \quad \fq_i\subseteq Q_i.$$
\end{enumerate}
With $Q_i,\, 1\leq i\leq r$,  as above,  $\Att_R(\lc^n_\fa(M))=T$
where $\fa=\bigcap\limits_{i=1}^rQ_i$.
\end{prop}

\begin{proof}
$(i)\Rightarrow (ii)$. By \cite[Theorem A]{DY1},
$\lc^n_\fa(R/\fp)\neq 0$ for all $\fp\in T$, that is $\fa+\fp$ is
$\fm$--primary for all $\fp\in T$ (by \LH). On the other hand, for
$1\leq i\leq r, \fq_i\notin T$ which is equivalent to say that
$\fa+\fq_i$ is not an $\fm$--primary ideal. Hence there exists a
prime ideal $Q_i\in \Supp_R(M)$ such that $\dim_R(R/Q_i)=1$ and
$\fa+\fq_i\subseteq Q_i$. It follows that $\underset{\fp\in
T}{\bigcap}\fp\nsubseteq Q_i$.
\\
\indent $(ii)\Rightarrow (i)$. Set $\fa:=\bigcap\limits_{i=1}^rQ_i$.
For each $i, 1\leq i\leq r$, $\fa+\fq_i\subseteq Q_i$ implies that
$\fa+\fq_i$ is not $\fm$--primary and so $\lc^n_\fa(R/\fq_i)= 0$.
Thus $\Att_R\lc_\fa ^n(M)\subseteq T$. Assume $\fp\in T$ and $Q\in
\Supp(M)$ such that $\fa+\fp \subseteq Q$. Then $Q_i\subseteq Q$ for
some $i, 1\leq i\leq r$. Since $\fp\nsubseteq Q_i$, we have $Q_i\neq
Q$, so $Q=\fm$. Hence $\fa+ \fp$ is $\fm$--primary ideal. Now, by
\LH, and by \cite[Theorem A]{DY1}, it follows that
$\fp\in\Att_R(\lc^n_\fa (M))$.
\end{proof}
\begin{cor}
If $\lc^n_\fa (M)\not=o$ then there is an ideal $\fb$ of $R$ such
that $\dim_R(R/\fb)\leq 1$ and $\lc^n_\fa (M)\cong\lc^n_\fb (M)$.
\end{cor}
\begin{proof}
 If $\Att_R(\lc^n_\fa(M))= \Assh_R(M)$, then $\lc^n_\fa (M)=
\lc^n_\fm (M)$. Otherwise $n\geq 1$ and $\Att_R(\lc^n_\fa(M))$ is a
proper subset of $\Assh_R(M)$. Set $\Assh_R(M)\setminus
\Att_R(\lc^n_\fa(M)):=\{\fq_1, \cdots, \fq_r\}$. By Proposition 2.1,
there are $Q_i\in \Supp_R(M)$ with $\dim_R(R/Q_i)= 1, \ i=1, \cdots,
r$,\ such that $\Att_R(\lc^n_\fa(M))= \Att_R(\lc^n_\fb(M))$ with
$\fb= {\bigcap\limits_{i=1}^rQ_i}$. Now, by \cite[Theorem 1.6]{DY2},
we have $\lc^n_\fa(M)\cong \lc^n_\fb(M)$. As $\dim(R/\fb)= 1$, the
proof is complete.
\end{proof}
\begin{cor}
If $\dim_R(M)=1$
 then any subset $T$ of $\Assh_R(M)$ is equal to the set
 $\Att_R(\lc^1_\fa(M))$ for some ideal $\fa$ of $R$.
 \end{cor}
 \begin{proof}
 With notations as in Proposition 2.1, we take $Q_i=\fq_i$ for $i=1,\cdots, r$.
  \end{proof}
  \indent By a straightforward argument one may notice that the condition ``complete" is
  superficial, for if $T$ is a non--empty proper subset of $\Assh_R(M)$,
  then $T=\Att_R(\lc^1_\fa(M))$, where
  $\fa=\underset{\fp\in\Assh_R(M)\setminus T}{\cap}\fp$.

 \vspace{0.5cm}
The following is an example to Proposition 2.1.
\begin{exam}
 Set $R=k[[X,Y,Z,W]]$, where $k$ is a field
and $X,Y,Z,W$ are
independent indeterminates. Then $R$ is a complete Noetherian
local ring with maximal ideal $\fm=(X,Y,Z,W)$. Consider prime
ideals
$$\fp_1=(X,Y) \quad  ,  \quad \fp_2=(Z,W)\quad , \quad  \fp_3=(Y,Z)
\quad , \quad  \fp_4=(X,W)$$ and  set  $\displaystyle
M=\frac{R}{\fp_1\fp_2\fp_3\fp_4}$ as an $R$--module, so that we have
$\Assh_R(M)=\{\fp_1,\fp_2,\fp_3,\fp_4\}$ and $\dim_R(M)=2$. We get
$\{\fp_i\}=\Att_R(\lc^2_{\fa_i}(M))$, where $\fa_1=\fp_2,
\fa_2=\fp_1, \fa_3=\fp_4, \fa_4=\fp_3$, and
$\{\fp_i,\fp_j\}=\Att_R(\lc^2_{\fa_{ij}}(M))$, where
$$\begin{array}{l}
\fa_{12}=(Y^2+YZ,Z^2+YZ,X^2+XW,W^2+WX),\\
\fa_{34}=(Z^2+ZW,X^2+YX,Y^2+YX,W^2+WZ),\\
\fa_{13}=(Z^2+XZ,W^2+WY,X^2+XZ),\\
\fa_{14}=(W^2+WY,Z^2+ZY,Y^2+YW),\\
\fa_{23}=(X^2+XZ,Y^2+WY,W^2+ZW),\\
\fa_{24}=(X^2+XZ,Y^2+WY,Z^2+ZW).\\
\end{array}$$
Finally, we have
$\{\fp_i,\fp_j,\fp_k\}=\Att_R(\lc^2_{\fa_{ijk}}(M))$, where
$\fa_{123}=(X,W,Y+Z)$, $\fa_{234}=(X,Y,W+Z)$, $\fa_{134}=(Z,W,Y+X)$.\\
\end{exam}

\begin{lem}
Assume that $n:=\dim_R(M)\geq 2$, and that $T$ is a non-empty subset
of $\Assh_R(M)$ such that
  $\underset{\fp\in T}{\bigcap} \fp \nsubseteq \underset
   {\fq\in \Assh_R(R/\sum\limits_{\fp\in T'}\fp)}{\bigcap} \fq$, where
   $T'=\Assh_R(M)\setminus T$.
Then there exists a prime ideal $Q\in \Supp_R(M)$ with
$\dim_R(R/Q)=1$ and $\Att_R(\lc^n_Q(M))=T.$
\end{lem}

\begin{proof}
Set $s:=\h_M(\sum\limits_{\fp\in T'}\fp)$. We have $s\leq n-1$,
otherwise $\Assh_R(R/\sum\limits_{\fp\in T'}\fp)= \{\fm\}$ which
contradicts the condition  $\underset{\fp\in T}{\bigcap} \fp
\nsubseteq \underset
   {\fq\in \Assh_R(R/\sum\limits_{\fp\in T'}\fp)}{\bigcap} \fq$.
As $R$ is catenary, we have $\dim_R(R/\sum\limits_{\fp\in
T'}\fp)=n-s$. We first prove, by induction on $j$, $0\leq j\leq
n-s-1$, that there exists a chain of prime ideals $Q_0 \subset Q_1
\subset \cdots \subset Q_j \subset \fm$ such that
$Q_0\in\Assh_R(R/\sum\limits_{\fp\in T'}\fp)$, $\dim_R(R/Q_j)=n-s-j$
and $\underset{\fp\in T}{\bigcap}\fp\nsubseteq Q_j$. There is
$Q_0\in\Assh_R(R/\sum\limits_{\fp\in T'}\fp)$ such that
$\underset{\fp\in T}{\bigcap}\fp\nsubseteq Q_0$. Note that
$\dim_R(R/Q_0)=\dim_R(R/\sum\limits_{\fp\in T'}\fp)=n-s$. Now,
assume that $0<j\leq n-s-1$ and that we have proved the existence of
a chain $Q_0 \subset Q_1 \subset \cdots \subset Q_{j-1}$ of prime
ideals such that $Q_0\in\Assh_R(R/\sum\limits_{\fp\in T'}\fp)$,
$\dim_R(R/Q_j)=n-s-(j-1)$ and that $\underset{\fp\in
T}{\bigcap}\fp\nsubseteq Q_{j-1}$. Note that we have
$n-s-(j-1)=n-s+1-j\geq 2$. Therefore the set $V$ defined as
\\
$$\begin{array}{ll}
V= \{\fq\in \Supp_R(M) |&  Q_{j-1}\subset \fq\subset \fq'\subseteq
 \fm, \dim_R(R/\fq)=n-s-j,\\ &
\fq'\in\Spec(R)\, \mbox{and}\,  \dim_R(R/\fq')=n-s-j-1\}
\end{array}$$
\\
is non-empty and so, by Ratliff's weak existence theorem
\cite[Theorem 31.2]{M}, is not finite. As $\underset{\fp\in
T}{\bigcap}\fp \nsubseteq Q_{j-1}$, we have $Q_{j-1}\subset
Q_{j-1}+\underset{\fp\in T}{\bigcap}\fp$. If, for $\fq\in V$,
$\underset{\fp\in T}{\bigcap}\fp\subseteq \fq$, then $\fq$ is a
minimal prime of $Q_{j-1}+\underset{\fp\in T}{\bigcap}\fp$. As $V$
is an infinite set, there is $Q_j\in V$ such that $\underset{\fp\in
T}{\bigcap}\fp\nsubseteq Q_j$. Thus the induction is complete. Now
by taking $Q:=Q_{n-s-1}$ and by Proposition 2.1, the claim follows.
\end{proof}

\begin{cor}
Assume that $n:=\dim_R(M)\geq 2$ and $T$ is a non-empty subset of
$\Assh_R(M)$ with $|T|=|\Assh_R(M)|-1$. Then there is an ideal $\fa$
of $R$ such that $\Att_R(\lc^n_\fa(M))=T$.
\end{cor}

\begin{proof}
Note that $\Assh_R(M)\setminus T$ is a singleton set $\{\fq\}$, say,
and so $\h_M(\fq)=0$ and $\underset{\fp\in T}{\bigcap}\fp\nsubseteq
\fq$. Therefore, by Lemma 2.5, the result follows.
\end{proof}

\begin{lem}
Assume that $n:=\dim_R(M)\geq 2$ and $\fa_1$ and $\fa_2$ are ideals
of $R$. Then there exists an ideal $\fb$ of $R$ such that
$\Att_R(\lc^n_\fb(M))=\Att_R(\lc^n_{\fa_{1}}(M))\cap\Att_R(\lc^n_{\fa_{2}}(M))$.
\end{lem}

\begin{proof}
Set $T_{1}=\Att_R(\lc^n_{\fa_{1}}(M))$ and
$T_{2}=\Att_R(\lc^n_{\fa_{2}}(M))$. We may assume that $T_1\bigcap
T_2$ is a non--empty proper subset of $\Assh_R(M)$. Assume that $\fq
\in \Assh_R(M)\setminus (T_1\bigcap T_2)=(\Assh_R(M)\setminus
T_1)\bigcup(\Assh_R(M)\setminus T_2) $. By Proposition 2.1, there
exists $Q\in \Supp_R(M)$ with $\dim_R(R/Q)=1$ such that
$\fq\subseteq Q$ and $\bigcap_{\fp\in T_1\bigcap T_2}\fp \nsubseteq
Q$. Now, by Proposition 2.1, again there exists an ideal $\fb$ of
$R$ such that $\Att_R(\lc^n_\fb(M))=T_1\bigcap T_2$.
\end{proof}

\indent Now we are ready to present our main result.

\begin{thm}
Assume that $T\subseteq \Assh_R(M)$, then there exists an ideal
$\fa$ of $R$ such that $T=\Att_R(\lc^n_\fa(M))$.
\end{thm}

\begin{proof}
By Corollary 2.3, we may assume that $\dim_R(M)\geq 2$ and that $T$
is a non-empty proper subset of $\Assh_R(M)$. Set
$T=\{\fp_1,\ldots,\fp_t\}$ and $\Assh_R(M)\setminus
T=\{\fp_{t+1},\ldots,\fp_{t+r}\}$. We use induction on $r$. For
$r=1$, Corollary 2.6 proves the first step of induction. Assume that
$r>1$ and that the case $r-1$ is proved. Set
$T_1=\{\fp_1,\ldots,\fp_t,\fp_{t+1}\}$ and
$T_2=\{\fp_1,\ldots,\fp_t,\fp_{t+2}\}$. By induction assumption
there exist ideals $\fa_1$ and $\fa_2$ of $R$ such that
$T_1=\Att_R(\lc^n_{\fa_1}(M))$ and $T_2=\Att_R(\lc^n_{\fa_2}(M))$.
Now by the Lemma 2.7 there exists an ideal $\fa$ of $R$ such that
$T=T_1\bigcap T_2=\Att_R(\lc^n_{\fa}(M))$.
\end{proof}
\begin{cor}
(See \cite[Corollary 1.7]{C}) With the notations as in Theorem 2.8,
the number of non-isomorphic top local cohomology modules of $M$
with respect to all ideals of $R$ is equal to $2^{|\Assh_R (M)|}$.
\end{cor}
\begin{proof}
It follows from Theorem 2.8 and \cite[Theorem 1.6]{DY2}.
\end{proof}
{\bf Acknowledgment.}The authors would like to thank the referee for
her/his comments.

\providecommand{\bysame}{\leavevmode\hbox
to3em{\hrulefill}\thinspace}

\end{document}